\documentclass[12pt]{amsart}

\newtheorem{Rmk}{Remark}[section]
\newtheorem{Theorem}[Rmk]{Theorem}

\newtheorem{Lemma}[Rmk]{Lemma}

\newtheorem{Prop}[Rmk]{Proposition}

\theoremstyle{plain}

\def\C{\mathbb C}

\def\R{\mathbb{R}}

\def\Z{\mathbb{Z}}

\def\to{\rightarrow}
\def\*{\times}

\begin{document}
\title[Finiteness of a section of the $SL(2,\C)$-character variety]
{Finiteness of a section of the $SL(2,\C)$-character variety 
of knot groups}
\author{Fumikazu Nagasato}
\address{Department of Mathematics, Tokyo Institute of Technology, 
2-12-1 Oh-okayama, Meguro-ku, Tokyo 152-8551, Japan}
\email{fukky@math.titech.ac.jp}

\begin{abstract}
We show that for any knot there exist only finitely many irreducible 
metabelian characters in the $SL(2,\C)$-character variety of the knot group, 
and the number is given explicitly by using the determinant of the knot. 
Then it turns out that for any 2-bridge knot 
a section of the $SL(2,\C)$-character variety 
consists entirely of all the metabelian characters, 
i.e., the irreducible metabelian characters and 
the single reducible (abelian) character. Moreover we find that 
the number of irreducible metabelian characters gives an upper bound 
of the maximal degree of the $A$-polynomial in terms of the variable $l$. 
\end{abstract}
\maketitle
\section{Introduction}\label{intro}
Our main object in this paper is representations of the knot group into $SL(2,\C)$. 
A representation $\rho$ of a group $G$ is said to be metabelian 
if the image $\rho([G,G])$ of the commutator subgroup $[G,G]$ 
is abelian. The irreducible metabelian representations are 
relatively easy to find in the irreducible representations in the following sense. 
Let us first denote the exterior of a knot $K$ in 3-sphere $S^3$ by $E_K$, 
which is the complement of an open tubular neighborhood of $K$ in $S^3$. 

\begin{Prop}\label{trace-free}
For an arbitrary knot $K$ in $S^3$, 
any irreducible metabelian representation $\rho$ of the knot group $\pi_1(E_K)$ 
into $SL(2,\C)$ satisfies 
\[ {\rm trace}(\rho(\mu))=0,\ {\rm trace}(\rho(\lambda))=2, \] 
where $(\mu,\lambda)$ is the pair of the elements in $\pi_1(E_K)$ 
represented by the standard meridian and longitude of the knot $K$, respectively. 
In particular, $\rho$ is not faithful. 
\end{Prop}
This proposition shows that all the irreducible metabelian characters 
of the knot group $\pi_1(E_K)$ are on the section of 
the $SL(2,\C)$-character variety cut by the equation 
$\chi_{\rho}(\mu)={\rm trace}(\rho(\mu))=0$, where $\chi_{\rho}$ is 
the character of $\rho$. 

In the paper \cite{l}, X.-S. Lin shows the finiteness of the number of 
conjugacy classes of the irreducible metabelian representations of 
the knot group into $SU(2)$. Moreover he gives the number of the conjugacy classes 
by using the Alexander polynomial $\Delta_K(t)$. 
We can in fact generalize the result to the case of $SL(2,\C)$. 

\begin{Theorem}\label{meta}
For any knot $K$ in $S^3$, 
there exist only finitely many conjugacy classes of the irreducible metabelian 
representations of the knot group $\pi_1(E_K)$ into $SL(2,\C)$. 
Moreover the number of the conjugacy classes (i.e., 
the number of irreducible metabelian characters) is given by 
\[ \frac{|\Delta_K(-1)|-1}{2}, \]
where 
the absolute value $|\Delta_K(-1)|$ 
is so-called the determinant of $K$. 
\end{Theorem}
For example, the number of irreducible metabelian characters is $\frac{p-1}{2}$ 
for the 2-bridge knot $S(p,q)$ (the Schubert presentation) 
because $|\Delta_{S(p,q)}(-1)|=p$ 
(see Proposition 8 in Gaebler's thesis \cite{g} for example). 
It turns out that the number of irreducible metabelian characters has 
an important meaning. 
For a finitely generated and finitely presented group $G$, let $X(G)$ be 
the $SL(2,\C)$-character variety of $G$ introduced in \cite{cs}. 
Then for a knot $K$ in $S^3$, we can define the $A$-polynomial 
$A_K(m,l)\in \Z[m,l]$ (refer to \cite{ccgls}). 
That is an invariant of knots constructed basically 
via the restriction 
\[ r:X(E_K):=X(\pi_1(E_K)) \to X(T^2):=X(\pi_1(T^2)) \] 
induced by the inclusion $i:\pi_1(T^2) \to \pi_1(E_K)$ 
(refer to \cite{ccgls} or Section $\ref{meaning}$). 

\begin{Theorem}\label{2-bridge}
For any 2-bridge knot $S(p,q)$, the number of irreducible metabelian 
characters in $X(\pi_1(E_{S(p,q)}))$ gives an upper bound 
of the maximal degree of the $A$-polynomial of $S(p,q)$ 
in terms of the variable $l$. In particular, the maximal degree is 
the number of irreducible metabelian characters surviving 
``discarding operations''. 
\end{Theorem}
We discuss the discarding operations 
at the very end of the proof of Theorem $\ref{2-bridge}$. 
Let $\deg_l(A_K(m,l))$ be the maximal degree of $A_K(m,l)$ in terms of $l$. 
Then the above theorem says that 
\[ \deg_l(A_{S(p,q)}(m,l))\leq \frac{p-1}{2}. \]
This gives a representation theoretical proof 
of the upper bound of $\deg_l(A_{S(p,q)}(m,l))$ shown in Le's paper \cite{le}. 

As a by-product, we get criteria of irreducible non-metabelian 
representations. For example, we have the following.

\begin{Prop}\label{A-poly}
If the $A$-polynomial of a knot $K$ satisfies 
$A_K(\pm\sqrt{-1},l)=0$, 
then there exist arcs in $X(E_K)$ each of which consists entirely of irreducible 
non-metabelian characters. 
If $A_K(m=\pm \sqrt{-1},l)\neq 0$ and has a factor other than $l-1$ or $l$ 
and moreover the knot $K$ is small, 
then there exists an irreducible non-metabelian representation satisfying 
${\rm trace}(\rho(\mu))=0$. 
\end{Prop}
We remark that the determinant of knots no longer controls the $A$-polynomial 
for knots other than 2-bridge knots. For example, see the knot $8_{20}$: 
\begin{eqnarray*}
&\deg_l(A_{8_{20}}(m,l))
=\deg_l(A_{8_{20}}(\pm \sqrt{-1},l))=\deg_l((l-1)^3(l+1)^2)=5,&\\
&\frac{|\Delta_{8_{20}}(-1)|-1}{2}=4.&
\end{eqnarray*}
(See also Section $\ref{remark}$.) 
In \cite{n}, we will introduce an algebraic variety
via the Kauffman bracket skein module. 
The number of its irreducible components 
will give upper bounds of both quantities $\frac{|\Delta_K(-1)|-1}{2}$ and 
$\deg_l(A_K(m,l))$ with some conditions. 

In this paper, we show the above statements in the following steps. 
In the next section, we give a proof of Theorem $\ref{meta}$, 
which includes a proof of Proposition $\ref{trace-free}$.  
In Section $\ref{meaning}$, we review the definition of the $A$-polynomial and 
give a proof of Theorem $\ref{2-bridge}$ and Proposition $\ref{A-poly}$. 
Then the finiteness of the section of the $SL(2,\C)$-character variety of 
the 2-bridge knot group is shown in the proof of Theorem $\ref{2-bridge}$. 
In the final section, we remark further aspects on this research. 
\section{Irreducible metabelian representations of knot groups}
We first consider a convenient presentation of knot groups 
to research the irreducible metabelian representations. 

\begin{Lemma}[Lemma 2.1 in \cite{l}]\label{pre}
The knot group $\pi_1(E_K)$ has a presentation
\[ \langle x_1,\cdots, x_{2g},\mu\ |\ \mu\alpha_i\mu^{-1}=\beta_i,i=1,\cdots,2g 
\rangle, \]
where $\mu$ is represented by a meridian of $K$, $\alpha_i$, $\beta_i$ are 
certain words in $x_1,\cdots,x_{2g}$, and $g$ is genus of a free Seifert surface 
of $K$. 
\end{Lemma}

This can be shown as follows (for more information refer to \cite{l}). 
Suppose $S$ be a free Seifert surface of a knot $K$. 
Namely, $(N(S),\overline{S^3-N(S)})$ is a Heegard splitting of $S^3$, 
where $N(S)$ is a closed tubular neighborhood of $S$ in $S^3$. 
Let $W_{2g}$ ($g={\rm genus}(S)$) be a spine of $S$, 
i.e., a bouquet (with a base point $*$) of circles in $S$ 
which is a deformation retract of $S$. 
Denote the oriented circles in $W_{2g}$ by $a_1,a_2,\cdots,a_{2g-1},a_{2g}$ 
forming a symplectic basis of $S$ with that order. Let $S\times [-1,1]$ be 
a bicollar of $S$ given by the positive normal direction of $S$ 
such that $S=S\times \{0\}$. Let $a_i^{\pm}=a_i\times \{\pm1\}$ 
be circles in $W_{2g}\times \{\pm 1\}$ respectively,  
for $i=1,\cdots,2g$. Note that there exists a circle $c$ on $\partial(S\times[-1,1])$, 
uniquely determined up to isotopy, 
such that $c$ connects two point $\{*\}\times\{\pm1\}$ and its interior 
has no intersection with $W_{2g} \times \{\pm1\}$. 
Consider the arc $c$ as a base point and 
choose a basis $x_1,\cdots,x_{2g}$ for a free group $\pi_1(S^3-S\times [-1,1])$. 
Then, for $i=1,\cdots,2g$, we define 
$\alpha_i$ (resp. $\beta_i$) by the element of $\pi_1(S^3-S\times[-1,1])$,
corresponding to $a_i^{+}$ (resp. $a_i^{-}$), 
which is a word in $x_1,\cdots,x_{2g}$.
The circle $c\cup\{*\}\times [-1,1]$ can be thought of as 
a meridian of a simple closed curve $K'$ on $S$ parallel to $K=\partial S$. 
Identifying $K'$ with $K$ and $\pi_1(S^3-S)$ with $\pi_1(S^3-S\times[-1,1])$, 
we can show Lemma $\ref{pre}$. 

In the above setting, 
we have the following fundamental facts on $\alpha_i$ and $\beta_i$, $i=1,\cdots,2g$, 
due to \cite{l}.  Let us denote by $v_{i,j}$ the exponent sum of $x_j$ in $\alpha_i$ 
and $u_{i,j}$ the exponent sum of $x_i$ in $\beta_i$. For two $2g\times 2g$ matrices 
$V:=(v_{i,j})$ and $U:=(u_{i,j})$ with integer entries, we have the followings:
\begin{itemize}
\item $V$ is so-called the Seifert matrix of the Seifert surface $S$,
\item $U=V^T$, where $T$ means transpose.
\end{itemize}
By the definition of the Alexander polynomial, we have 
\begin{eqnarray}\label{alex}
|\det(V+V^T)|=|\Delta_K(-1)|.
\end{eqnarray}
\noindent
{\bf Proof of Proposition $\ref{trace-free}$ and Theorem $\ref{meta}$}: 
The proof is straightforward. 
It is well-known that there exist only two maximal abelian subgroups 
Hyp(=hyperbolic) and Para(=parabolic) in 
$SL(2,\C)$, up to conjugation: 
\[{\rm Hyp}:=\left\{\left(\left.
\begin{array}{cc} 
\lambda & 0\\
0 & \lambda^{-1}
\end{array}\right)
\in SL(2,\C)\right|\lambda \in \C^{*}\right\} \]
\[ {\rm Para}:=\left\{\pm
\left(\left.
\begin{array}{cc}
1 & \omega\\
0 & 1
\end{array}\right)
\in SL(2,\C)\right|\omega \in \C\right\}
\]
Note that two non-trivial elements $g$ and $h$ in $SL(2,\C)$ 
are commutative if and only if 
the subgroup in $SL(2,\C)$ generated by $g$ and $h$ is conjugate to a subgroup 
of Hyp or Para. 
It is easy to see that the elements 
$x_1,\cdots,x_{2g}$ of $\pi_1(E_K)$ are in the commutator subgroup 
$[\pi_1(E_K),\pi_1(E_K)]$. 
Let $\rho$ be an arbitrary irreducible metabelian representation. 
Then we can assume up to conjugation that 
\[ (\rho(x_i))_{1\leq i \leq 2g}=\left(\left(
\begin{array}{cc}
\lambda_i & 0\\
0 & \lambda_i^{-1}
\end{array}\right)\right)_{1\leq i\leq 2g}
\mbox{ or }
\pm \left(\left(
\begin{array}{cc}
1 & \omega_i\\
0 & 1
\end{array}\right)\right)_{1\leq i \leq 2g}, \]
where $\lambda_i(\neq 0)$ and $\omega_i$ are complex numbers. 
We can also check that $\mu x_i\mu^{-1}$, $i=1,\cdots,2g$, are in the commutator 
subgroup. With some linear algebra, we see that $\rho$ is always an abelian 
representation if $\rho(x_i)$ are in Para. 
Hence we can ignore the parabolic case. In the hyperbolic case, we have
\[ \rho(\mu)=\left(
\begin{array}{cc}
0 & b\\
-b^{-1} & 0
\end{array}
\right),\ 
\rho(\mu x_i\mu^{-1})=\left(
\begin{array}{cc}
\lambda_i^{-1} & 0\\
0 & \lambda_i
\end{array}
\right),\mbox{ for $i=1,\cdots,2g$, }
\]
where $b$ is a non-zero complex number. Therefore, any irreducible metabelian 
representation satisfies ${\rm trace}(\rho(\mu))=0$. Note that 
the longitude $\lambda$ can be described as a word in $x_1,\cdots,x_{2g}$ 
with zero exponent sum. 
Since $x_i,\ i=1,\cdots,2g$, are in the commutator subgroup, we have 
\[ \rho(\lambda)=\left(\begin{array}{cc}1&0\\0&1\end{array}\right), \] 
and thus ${\rm trace}(\rho(\lambda))=2$. In particular, $\rho$ is not faithful.
This gives a proof of Proposition $\ref{trace-free}$. 

Now, substituting the above equations to the defining relations of $\pi_1(E_K)$, 
we have the following equations for $\lambda_1,\cdots,\lambda_{2g}$:
\begin{eqnarray}\label{eq}
\lambda_1^{w_{i,1}}\cdots\lambda_{2g}^{w_{i,2g}}=1, 
\mbox{ for $i=1,\cdots,2g$}, 
\end{eqnarray}
where $(w_{i,j})=V+V^T$, the absolute value $|\det(w_{i,j})|$ 
is the determinant $|\Delta_K(-1)|$ of $K$ (see Equation $(\ref{alex})$).  
Put $\lambda_i:=r_ie^{\sqrt{-1}\theta_i}$, $r_i>0$, $0\leq\theta_i<2\pi$, 
for $i=1,\cdots,2g$, to solve Equations $(\ref{eq})$. 
Then we get 
\begin{eqnarray}
r_1^{w_{i,1}}\cdots r_{2g}^{w_{i,2g}}=1\label{r},\\
e^{\sqrt{-1}(w_{i,1}\theta_1+\cdots+w_{i,2g}\theta_{2g})}=1.
\label{lambda}
\end{eqnarray}
We first solve Equations $(\ref{r})$. 
Taking a function $\log=\log_e$ for its both sides, we have
\[ (w_{i,j})(\log r_1,\cdots,\log r_{2g})^T=(0,\cdots,0)^T, \]
where $T$ is the transpose. 
Let us take the inverse matrix of $(w_{i,j})$ for both sides. 
Indeed, this can be successfully done because $|\det(w_{i,j})|=|\Delta_K(-1)|\neq 0$. 
Since the function $\log$ is bijective on the real numbers $\R$, 
we have only a single solution
\[ (r_1,\cdots,r_{2g})=(1,\cdots,1). \]
We next solve Equations $(\ref{lambda})$. They are equivalent to the following. 
Let us think of $\{\theta_i\}_{i=1,\cdots,2g}$ 
as elements $\{[\theta_i]\}_{i=1,\cdots,2g}$ of the quotient space 
$\R/(2\pi n \sim 2\pi m)\cong[0,2\pi)$. Then we have  
\[ w_{i,1}[\theta_1]+\cdots+w_{i,2g}[\theta_{2g}]=[0], 
\mbox{ for $i=1,\cdots,2g$}. \]
They can be encoded into a single form by using the matrix $(w_{i,j})$ as follows:
\begin{eqnarray}\label{overZ}
(w_{i,j})([\theta_1],\cdots,[\theta_{2g}])^T=([0],\cdots,[0])^T.
\end{eqnarray}
Then the fundamental fact on the order ideal gives the number of solutions 
of the above equations as the absolute value of the determinant of 
the matrix $(w_{i,j})$ (for example, see p.205 of \cite{r}). 
Recalling $|\det(w_{i,j})|=|\Delta_K(-1)|$, we see that 
$|\Delta_K(-1)|$ is the number of solutions of Equations $(\ref{lambda})$ 
and thus Equation $(\ref{eq})$. 

Now, we discard the trivial solution $(\lambda_1,\cdots,\lambda_{2g})=(1,\cdots,1)$ 
from $|\Delta_K(-1)|$ solutions since it gives an abelian representation. 
We remark that a non-trivial solution of 
Equations $(\ref{eq})$ always gives an irreducible representation, because 
a representation $\rho$ of a finitely generated and finitely presented group $G$ 
into $SL(2,\C)$ 
is reducible if and only if ${\rm trace}(\rho(g))=2$ for any $g\in[G,G]$ 
(see Corollary 1.2.1 of \cite{cs}). 
Remember that there still exist infinitely many possibility for $\rho$ to be 
an irreducible metabelian representation because of $b$ in $\rho(\mu)$. 
However they are all conjugate. Namely, we first fix a non-trivial solution 
$(\lambda_i)$. 
Let $\rho_{(\lambda_i),b}$ and $\rho_{(\lambda_i),b'}$ be representations 
associated with non-zero complex numbers $b$ and $b'$ for $\rho(\mu)$ 
along with the above fixed solution $(\lambda_i)$, respectively. 
Then we can easily show with some linear algebra that they are always conjugate. 
Hence we do not care about $b$ for the conjugacy classes. 
On the other hand, if $(\lambda_1,\cdots,\lambda_{2g})$ is a solution of Equations 
$(\ref{lambda})$, then $(\lambda_1^{-1},\cdots,\lambda_{2g}^{-1})$ is also 
a solution. That is, let $\rho_{(\lambda_i),b}$ and $\rho_{(\lambda_i'),b}$ be 
representations defined in the same fashion as above. Then it can be easily shown 
that they are conjugate if and only if $(\lambda_i')=(\lambda_i^{-1})$ or ($\lambda_i$).
Combining those facts on conjugation, we see that the number of conjugacy classes 
of the irreducible metabelian representations is 
$\frac{|\Delta_K(-1)|-1}{2}$.
$\hfill \square$
\section{Representation theoretical meaning of the maximal degree of 
the $A$-polynomial in terms of $l$}\label{meaning}
We first review the $A$-polynomial of knots (for more information, refer to 
\cite{ccgls}). For a finitely generated and finitely presented group $G$, 
the $SL(2,\C)$-character variety $X(G)$ is defined 
as the quotient of ${\rm Hom}(G,SL(2,G))$ 
by the trace function $t_g$, $g\in G$,
\[ t_g:{\rm Hom}(G,SL(2,\C))\to\C,\ t_g(\rho):={\rm trace}(\rho(g)), \]
i.e., two representations $\rho_1$ and $\rho_2$ are equivalent if 
$t_g(\rho_1)=t_g(\rho_2)$ for any $g\in G$. 
Actually, it turns out that the set $X(G)$ can be identified 
with an algebraic variety in some space $\C^N$. 

Now, for a knot $K$ in $S^3$, we consider the restriction 
\[ r:X(E_K):=X(\pi_1(E_K)) \to X(\partial E_K=T^2):=X(\pi_1(T^2)),
\ r(\chi_{\rho}):=\chi_{\rho\circ i} \] 
induced by the inclusion $i:\pi_1(T^2) \to \pi_1(E_K)$. 
Every irreducible component of $X(E_K)$ has 1- or 0-dimensional closure 
of the image under $r$ (refer to \cite{ccgls} and also Lemma 2.1 in \cite{dg}). 
Let $\Delta$ be the set of pairs of diagonal matrices 
\[ \Delta:=\left\{\left. \left( 
\left(\begin{array}{cc}m & 0\\ 0 & m^{-1}\end{array}\right),
\left(\begin{array}{cc}l & 0\\ 0 & l^{-1}\end{array}\right)
\right) \right| (m,l)\in\C^*\times\C^* \right\} \]
Then take a preimage $\mathcal{E}(K)$ of $r(X(E_K))$ under a 2-fold branched 
covering map $p:\Delta \to X(T^2)$ defined by 
\[ p\left(\left(\begin{array}{cc}m & 0\\ 0 & m^{-1}\end{array}\right),
\left(\begin{array}{cc}l & 0\\ 0 & l^{-1}\end{array}\right)\right)
:=(m+m^{-1},l+l^{-1},ml+m^{-1}l^{-1})\in\C^3. \]
Identifying $\Delta$ with the set $\C^*\times \C^*$ 
via the correspondence
\[ \left( 
\left(\begin{array}{cc}m & 0\\ 0 & m^{-1}\end{array}\right),
\left(\begin{array}{cc}l & 0\\ 0 & l^{-1}\end{array}\right)
\right) \to (m,l)\in\C^*\times\C^*, \] 
we can think of $\mathcal{E}(K)$ as a subset 
(also denoted by $\mathcal{E}(K)$) in $\C^*\times \C^*$. 
Taking the closure of $\mathcal{E}(K)$ in $\C^2$, we get an algebraic variety 
$\overline{\mathcal{E}(K)}$ in $\C^2$ (called the eigenvalue variety of knot $K$) 
each of whose irreducible components is 1- or 0-dimensional. 
Discard all the 0-dimensional components and denote the remains by $D_K$. 
The defining polynomial of $D_K$ in $\Z[m,l]$, uniquely determined up to 
non-zero constant multiple, is called the $A$-polynomial of knot $K$ 
and denoted by $A_K(m,l)$. 
Since $A_K(m,l)$ has always the factor $l-1$, which comes from the characters of 
abelian representations, we divide $A_K(m,l)$ by $l-1$ and again denote 
the resulting polynomial by $A_K(m,l)$. Note that the variable $m$ of $A_K(m,l)$ 
has only even power. 

The set $\overline{\mathcal{E}(K)}-\mathcal{E}(K)$ may have 
finitely many points. Such a point, say $(m_0,l_0)$, 
is called a hole of the eigenvalue variety $\overline{\mathcal{E}(K)}$ if 
that satisfies $m_0\neq 0$ and $l_0\neq 0$. 
We remark that if $\overline{\mathcal{E}(K)}$ has a hole, then the knot exterior 
has a closed essential surface (see Section 5 in \cite{c-l} for example). 

\noindent
{\bf Proof of Theorem $\ref{2-bridge}$}: 
Let us first recall the 2-bridge knots. For coprime odd integers $p$ and $q$, 
($p>0$, $p>|q|>0$), there exists associated the 2-bridge knot $S(p,q)$ 
(refer to \cite{k} for example). The knot group 
$\Gamma_{p,q}:=\pi_1(E_{S(p,q)})$ has the form 
\[ \Gamma_{p,q}=\langle x_1,x_2\ |\ wx_1=x_2w,\ 
w:=x_1^{e_1}x_2^{e_2}\cdots x_2^{e_{p-1}}\rangle, \]
where $e_i:=(-1)^{[\frac{iq}{p}]}$ for $0 \leq i \leq p-1$, $[\cdot]$ is 
the Gaussian integer. 
By the definition of the $A$-polynomial, 
it suffices to focus on the non-abelian representations to prove the theorem. 
It follows from Lemma 1 in \cite{ri} that 
any non-abelian representation $\rho:\Gamma_{p,q} \to SL(2,\C)$ can be conjugated 
so that 
\[ \rho(x_1)=
\left(\begin{array}{cc}t^{1/2}& t^{-1/2}\\0&t^{-1/2}\end{array}\right),\ 
\rho(x_2)=
\left(\begin{array}{cc}t^{1/2}& 0\\ -t^{1/2}u&t^{-1/2}\end{array}\right). 
\]
Put 
\[ \rho(w):=\left(\begin{array}{cc}w_{11}(t,u)&w_{12}(t,u)\\
w_{21}(t,u)&w_{22}(t,u)\end{array}\right),\ 
w_{ij}(t,u)\in\Z[t,t^{-1},u]. \]
Then Theorem 1 in \cite{ri} combined with some easy calculation shows that 
$\rho$ is a representation of $\Gamma_{p,q}$ into $SL(2,\C)$ if and only if 
$t$ and $u$ satisfy 
\begin{eqnarray}\label{rep}
w_{11}(t,u)+(1-t)w_{12}(t,u)=0. 
\end{eqnarray}
We now focus on the section $S_0(S(p,q))$ of $X(E_{S(p,q)})$ using the equation 
$\chi_{\rho}(\mu)={\rm trace}(\rho(\mu))=0$, where $\chi_{\rho}$ is 
the character of $\rho$. This means assuming $t=-1$, and thus 
\[ \rho(x_1)=
\left(\begin{array}{cc}\sqrt{-1}& -\sqrt{-1}\\0&-\sqrt{-1}\end{array}\right),\ 
\rho(x_2)=
\left(\begin{array}{cc}\sqrt{-1}& 0\\ -\sqrt{-1}u&-\sqrt{-1}\end{array}\right). 
\]
Then Equation ($\ref{rep}$) goes to 
\begin{eqnarray}\label{rept}
w_{11}(-1,u)+2w_{12}(-1,u)=0.
\end{eqnarray}
Since $\rho(x_i)^{-1}=-\rho(x_i)$ at $t=-1$ for $i=1,2$, 
and $e_j=e_{(p-1)-j}$ for $1 \leq j \leq p-1$, we have 
\[ \rho(w)=(\rho(x_1)\rho(x_2))^{\frac{p-1}{2}}=
\left(\begin{array}{cc}-1-u &-1\\-u&-1\end{array}\right)^{\frac{p-1}{2}}. 
\]
By induction, we get $\deg_u(w_{11}(-1,u))=\frac{p-1}{2}$, 
$\deg_u((w_{12})(-1,u))=\frac{p-3}{2}$. Hence Equation $(\ref{rept})$ have 
$\frac{p-1}{2}$ solutions over $\C$ with multiplicity. 
Remember that there exist $\frac{|\Delta_{S(p,q)}(-1)|-1}{2}=\frac{p-1}{2}$ 
characters of 
irreducible metabelian representations on the section $S_0(S(p,q))$. 
So all $\frac{p-1}{2}$ solutions of Equation ($\ref{rept}$) must be distinct and 
associated representations are irreducible metabelian and 
become $\frac{p-1}{2}$ distinct points on $S_0(S(p,q))$. 
Namely, $S_0(S(p,q))$ consists of the abelian character and 
$\frac{p-1}{2}$ points corresponding to the irreducible metabelian characters. 
Note that for any knot $K$, $S_0(K)$ has only irreducible characters 
except the abelian character, since $\Delta_K(-1)\neq 0$ (refer to \cite{hpp} 
or the original papers \cite{b,d}). 

Now, for any 2-bridge knot the eigenvalue variety has no holes, 
because any 2-bridge knot is small, namely the knot exterior 
has no closed essential surfaces. 
Moreover Proposition 2 in the paper of Hatcher and Thurston \cite{ht} shows 
that $S(p,q)$ has no meridional boundary slope. 
Then it follows from Theorem 3.4 in \cite{ccgls} 
that there exist no vertical edges in the Newton polygon of 
the $A$-polynomial of $S(p,q)$. Hence $A_{S(p,q)}(\sqrt{-1},l)$ must have 
constant term and thus $(m,l)=(\sqrt{-1},0)$ is not a solution of 
$A_{S(p,q)}(\sqrt{-1},l)=0$. Therefore every solution of 
$A_{S(p,q)}(\sqrt{-1},l)=0$ comes from $S_0(S(p,q))$. 
This means that the number of solutions of $A_{S(p,q)}(\sqrt{-1},l)=0$ 
with multiplicity, which is equal to the maximal degree of $A_{S(p,q)}(\sqrt{-1},l)$, 
must be less than or equal to the number of irreducible metabelian characters 
in $S_0(S(p,q))$. Note that the above inequality can be caused by 
the fact that we do not consider the 0-dimensional components of $D_K$ 
and we define the $A$-polynomial so that the polynomial has no repeated factors. 

Now recall that the presentation of longitude $\lambda$ using the generators 
$x_1$, $x_2$ is 
\[ \lambda=\omega^{-1} \cdot \widetilde{\omega} \cdot x_1^{2\sigma}, \]
where $\widetilde{\omega}:=x_1^{-\varepsilon_1}x_2^{-\varepsilon_2}\cdots
x_2^{-\varepsilon_{p-1}}$, 
$\sigma:=\sum_{i=1}^{p-1}\varepsilon_i$.
Since $\rho(x_i)^{-1}=-\rho(x_i)$ for $i=1,2$, we have $\rho(\lambda)={\rm id}$. 
So the $A$-polynomial at $m=\sqrt{-1}$ has the factorization 
\[ A_{S(p,q)}(\sqrt{-1},l)=(l-1)^{k_{p,q}},\ k_{p,q}\in\Z_{\geq 0}. \] 
This shows that 
\[ \frac{p-1}{2}\geq k_{p,q}=\deg_l(A_{S(p,q)}(\sqrt{-1},l)). \]
Since there exist no vertical edges in the Newton polygon of $A_{S(p,q)}(m,l)$, 
we have 
\[ \deg_l(A_{S(p,q)}(\sqrt{-1},l))=\deg_l(A_{S(p,q)}(m,l)). \]
Therefore we get 
\[ \frac{p-1}{2} \geq k_{p,q}=\deg_l(A_{S(p,q)}(\sqrt{-1},l))
=\deg_l(A_{S(p,q)}(m,l)). \]
This shows the first statement in Theorem $\ref{2-bridge}$.

As regards the meaning of the maximal degree of the $A$-polynomial of $S(p,q)$, 
the process of throwing away the 0-dimensional components of $D_K$ 
can make a difference between $\frac{p-1}{2}$ and $\deg_l(A_{S(p,q)}(m,l))$. 
We also remark that the $A$-polynomial is defined so that 
there are no repeated factors in it, 
namely the multiplicity of the curve components in $D_K$
are ignored in the construction process of the $A$-polynomial. 
This discarding process can also decrease the number of irreducible 
metabelian characters in the target. 
So we have to count the number of irreducible metabelian characters 
surviving the above {\it discarding operations}
when calculating the exact maximal degree of $A_{S(p,q)}(m,l)$ 
in terms of $l$. 
Hence we can say that $\deg_l(A_{S(p,q)}(m,l))$ is 
the number of irreducible metabelian 
characters surviving the discarding operations. 
$\hfill \square$

\noindent
{\bf Proof of Proposition $\ref{A-poly}$}: 
Let us first factorize $A_K(\pm \sqrt{-1},l)$ over $\C$. 
If 
\[ A_K(\pm \sqrt{-1},l)=0, \] 
then $A_K(m,l)$ must have a factor $m^2+1$, since 
$m^2+1$ is the minimal polynomial of $\sqrt{-1}$ over $\Z$. 
Then there exist two irreducible components $m=\pm \sqrt{-1}$ of $D_K$ in $\C^2$, 
denoted respectively by $C_{\pm}$. By Proposition $\ref{meta}$, 
all the points of $C_{\pm}$ except $(m,l)=(\pm \sqrt{-1},1)$ and holes 
correspond to 
characters of irreducible non-metabelian representations. Hence there must exist 
arcs in $X(E_K)$ associated with $C_{\pm}-\{(\pm\sqrt{-1},1)\}\cup\{\mbox{holes}\}$ 
consisting entirely of characters of irreducible non-metabelian representations. 
This completes the proof of the first criterion.

As regarding the second criterion, 
we assume that there exists a factor $l-\omega$, $\omega(\neq 0,1) \in \C^*$. 
Note that the points $(m,l)=(\pm \sqrt{-1},\omega)$ are not holes since 
$K$ is small. 
It follows from Proposition $\ref{trace-free}$ that 
the characters of all the irreducible metabelian representations surviving 
the discarding operations, as well as all the abelian characters, 
correspond to the factor $l-1$. 
Hence we see that there exists at least a character coming from 
an irreducible non-metabelian representation $\rho$ surviving the discarding 
operations such that ${\rm trace}(\rho(\mu))=0$, 
${\rm trace}(\rho(\lambda))=\omega+\omega^{-1}$. 
This completes the proof of the second criterion. 

\hfill $\square$
\section{Remarks}\label{remark}
As remarked in Section $\ref{intro}$, the determinant of knots no longer controls 
the $A$-polynomial for knots other than 2-bridge knots. 
However it is possible that the quantity $\frac{|\Delta_K(-1)|-1}{2}$ 
will be able to give an upper bound of the multiplicity of the factor 
$l-1$ of $A_K(\pm \sqrt{-1},l)$ 
for any $K$. To check this, one may want to observe whether or not 
every character $\chi_{\rho}$ on the section $S_0(K)$ satisfying 
$\chi_{\rho}(\lambda)=2$ is a metabelian character. 

As performed in Section $\ref{meaning}$, the fact that the Newton polygon of 
$A_{S(p,q)}(m,l)$ has no vertical edges as well as the fact that 
there exist no holes in $\overline{\mathcal{E}(S(p,q))}$ is very useful to study 
the $A$-polynomial by using the section $S_0(S(p,q))$. 
It will be very interesting to research 
whether or not there exist vertical edges in the Newton polygon of $A_K(m,l)$ 
and there exists a hole in $\overline{\mathcal{E}(K)}$ 
for knots other than 2-bridge knots. 

In this paper, we have a little discussion on the exact degree 
of the $A$-polynomial of the 2-bridge knots by using the word 
``discarding operations''. 
We will discuss more this topic at another time. 
\section*{Acknowledgement}
The author would like to thank Professor Xiao-Song Lin for very fruitful 
discussions and useful comments. 
He is also grateful to Masaharu Ishikawa for helpful discussions. 
Finally, the author has been supported by JSPS Research Fellowships 
for Young Scientists.

\end{document}